\documentclass{amsart}
\usepackage{graphicx}
\theoremstyle{definition}

\theoremstyle{remark}

\numberwithin{equation}{section}

\begin {document}
\title{Discrete Algebraic Equations and Discrete Operator Equations(Presentations for ICM 2010)}
\author{Wu Zi qian}
\address{Fangda group company,Shenzhen city,Guangdong province,China}
\email{$runton_-runton$$@$ruc.edu.cn,woodschain$@$sohu.com}  \pagestyle{myheadings}
\begin{abstract}
\mbox{} We give constructive results of Hilbert's 13th problem for
discrete functions. By them we give formula solution expressed by a
superposition of functions of one variable to equations constructed
by discrete functions and equations with parameterized discrete
functions. Further more we give formula solution expressed by a
superposition of operators of one variable to equations constructed
by discrete operators and equations with parameterized discrete
operators. This is a Short communication, Section 9,Functional Analysis and Application, Saturday, August 21, 2010,18:00-18:15, Room No. T3.
\end{abstract}
\maketitle
\section{Introduction}

Problems about equations are very important and difficult. Solving
quadratic equation and cubic equation and quartic equation had cost
the mathematicians in history a great deal of time.

Babylonians solve quadratics in radicals in 2000 BC. Cubic equation
and quartic equation were solved by Italian mathematicians
Girolamo.Cardano(1501-1576) and
Ludovico.Ferrari(1522-1565)in 16th century,respectively.

But mathematicians met big troubles when they tried to solve quintic
equation. Leonhard.Euler(1707-783) believed quintic equation
can be changed to a quartic equation by transformation of variable.
Niels.Henrik.Abel (1802-1829) got a conclusion that there is no
solution by radicals for a general polynomial algebraic equation if
n$\geq$5. Evariste.Galois (1811-1832) built group theory and
got the same conclusion. His method come down to now and can be
found in any textbook about Galois group theory.

There is no solution by radicals. Are there any solutions of other
forms such as numerical solution and solution expressed in function
of two variables or of many variables or solution expressed in
series or in integral?

We do not discuss numerical solutions because they belong to applied
mathematics. We prefer formula solution expressed in binary function
to other ones. What is a formula solution expressed in binary
function? It contains only function of two variables. We can give a
expression of a alone binary function at the beginning. We can
replace any one of variables by a binary function then we get a new
expression. We can replace any one of variables of this new
expression by a binary function again and get a more complex
expression. We can repeat the procedure for any finite times. But it
is not easy to get solution expressed in binary function.  It's
easier to get solutions of other forms.History developed just like
this.

Camille.Jordan (1838-1922) shows that algebraic equations of
any degree can be solved in terms of modular functions in 1870.
Ferdinand.von.Lindemann (1852-\\1939) expresses the roots of an
arbitrary polynomial in terms of theta functions in 1892. In 1895
Emory.McClintock (1840-1916) gives series solutions for all the
roots of a polynomial. Robert Hjalmal.Mellin (1854-1933) solves
an arbitrary polynomial equation with Mellin integrals in 1915. In
1925 R.Birkeland shows that the roots of an algebraic equation can
be expressed using hypergeometric functions in several variables.
Hiroshi.Umemura expresses the roots of an arbitrary polynomial
through elliptic Siegel functions in1984[1].

All of solutions mentioned above are not ones expressed in binary
function.

By Tschirnhausen transformation a quintic equation or a sextic
equation can be changed to ones containing only two parameters so
there are solutions expressed in binary function for them.
David.Hilbert presumed that there is no solution expressed in binary
function for polynomial equations of n when n$\geq$7 and wrote his
doubt into his famous 23 problems as the 13th one[2].

Hilbert published his last mathematical paper [3] in 1927 where he
reported on the status of his problems, he devoted 5 pages to the
13th problem and only 3 pages to the remaining 22 problems. We can
see that so much attention Hilbert paid to 13th problem. In 1957
V.I.Arnol'd proved that every continuing function of many variables
can be represented as a superposition of functions of two variables
and refuted Hilbert conjecture[4][5]. Furthermore, A.N.Kolmogorov
proved that every continuous function of several variables can be
represented as a superposition of continuous functions of one
variable and the operation of addition [6].

Result for Hilbert's 13th problem is very important for us and it
points us a quite right direction to solve polynomial equations and
general algebraic equations. But method used in it is topological and
the result is not a constructive one. In this paper we will give a
constructive result in discrete situation. This result is very
important. We can construct profuse discrete algebraic equations and
discrete operator equations and for this result we can give any of
them a formula solution. There is never such a mathematical
structure in the history of mathematics. This is the first time!
A.G.Vitushkin dissatisfies the current results about 13th problem
and points out that the algebraic core of the problem remains
untouched[7]. We believe we have gotten the algebraic core Vitushkin
wanted.

\section{Constructive results for
Hilbert's 13th problem}

A.N.Kolmogorov expresses function of several variables as a
superposition of functions of one variable like this:
\begin {align}\label{eq:eps}
W(x_{1},x_{2}\cdots,x_{n})=\sum_{i=1}^{2n+1}f_{i}\Big[g_{i1}(x_{1})+
g_{i2}(x_{2})+\cdots+g_{in}(x_{n})\Big]
\end{align}
This is a existence  result but it's easy to give a constructive
result for discrete functions.

{\bf Definition~2.1}
Let A=\{-1,0,1\}, a three numbers function of M
variables is defined as:

\qquad \qquad \qquad \qquad\qquad g:$A^{M}$$\longrightarrow$A

There are 9 discrete points for a binary three numbers function. A
binary three numbers  function can be indicated by a table with 4x4
elements. Its first column indicates the first variable and the
first row indicates the second variable. To give a table with 4x4
elements is to define a binary three numbers  function and vice
versa. For example:

 \begin{tabular}{|c|c|c|c|}
  \hline
       &-1& 0& 1 \\
  \hline
       -1&1&-1&0\\
  \hline
       0&-1&0&1\\
  \hline
       1&0&1&-1\\
  \hline
\end{tabular}
{          }
  \qquad \qquad\begin{tabular}{|c|c|c|c|}
  \hline
       &-1& 0& 1 \\
  \hline
       -1&0&0&0\\
  \hline
       0&0&0&0\\
  \hline
       1&0&0&0\\
  \hline
\end{tabular}
{          }
 \qquad \qquad\begin{tabular}{|c|c|c|c|}
  \hline
       &-1& 0& 1 \\
  \hline
       -1&1&1&1\\
  \hline
       0&0&0&0\\
  \hline
       1&-1&-1&-1\\
  \hline
\end{tabular}

There are three functions in above tables. The first one is linear
binary three numbers  function and the second one is identity
function with o value and value of the third one is not change with
the second variable. There is only one value for each discrete point
in these three functions and they are called single-valued binary
three numbers function. It's easy to know there are $3^{9}$=19683
single-valued binary three numbers  function. Can it be two-valued
or three-valued in a discrete point?certainly! There are three
combinations -1,0 and -1,1 and 0,1 for two-valued and only one
combination-1,0,1 for three-valued and numbers will be partitioned
by symbol '*' if it's a multi-valued. Can it be no-valued in a
discrete point? Yes! We will indicate it in 'N'. A binary three
numbers function can be no-valued in all 9 discrete points in
uttermost like:

\qquad \qquad\qquad \qquad\qquad\qquad
\begin{tabular}{|c|c|c|c|}
  \hline
       &-1& 0& 1 \\
  \hline
       -1&N&N&N\\
  \hline
       0&N&N&N\\
  \hline
       1&N&N&N\\
  \hline
\end{tabular}

There is single-valued,two-valued,three-valued and no valued point
in below three numbers function.

\qquad\qquad \qquad\qquad \qquad
\begin{tabular}{|c|c|c|c|}
  \hline
       &-1& 0& 1 \\
  \hline
       -1&-1&0&1\\
  \hline
       0&-1*0&-1*1&0*1\\
  \hline
       1&N&-1*0*1&N\\
  \hline
\end{tabular}

It's easy to know there are $8^{9}$ binary three numbers  functions.
A unary three numbers function  can be indicated by three value
numbers partitioned by symbol ',' in bracket and numbers will be
partitioned by the symbol '*'if it's many-valued, for
example:(-1*0,N,-1*0*1).

The expression $H(x_{1},x_{2})=f[g_{11}(x_{1})+g_{12}(x_{2})]$
contains $x_{1}$,$x_{2}$, but we intend to take H as a independent
object not containing  $x_{1}$,$x_{2}$. We  can't express H in
$f[g_{11}+g_{12}]$ because we will get a unary function
$[g_{11}+g_{12}]$ by adding $g_{11}$ and $g_{12}$.
$f[g_{11}+g_{12}]$ is also a unary function and is never equal to
binary the function H.

We express H by only f, $g_{11}$, $g_{12}$ without $x_{1},
x_{2}$like this:

\qquad \qquad \qquad \qquad \qquad
$H=f[g_{11}\widetilde{\alpha_{1}}+g_{12}\widetilde{\alpha_{2}}]$

To define a function is to give a rule to get it's values. For such
an expression we are very clear the rule about getting values of the
function if we replace $\widetilde{\alpha_{1}}$ or
$\widetilde{\alpha_{2}}  $by $x_{1} or x_{2}$ ,respectively. That is
enough.

Binary three numbers  function is called single term binary three
numbers function if it can be represented as
$H=f[g_{11}\widetilde{\alpha_{1}}+g_{12}\widetilde{\alpha_{2}}]$ in
which f ,$g_{ij}$ is unary three numbers function and it will be
called L term binary three numbers  function if it can be expressed
as $\sum
f_{i}[g_{i1}\widetilde{\alpha_{1}}+g_{i2}\widetilde{\alpha_{2}}]
(i=1,L)$. Expressing a function of many variables as this form is
also called representing it as a superposition of functions of one
variable or decomposing it to functions of one variable.

For example F is a single term binary three numbers  function:

\qquad\qquad \qquad \qquad
$F=(0,0,1)\Big[(0,0,1)\widetilde{\alpha_{1}}+(-1,-1,0)\widetilde{\alpha_{2}}\Big]$

{\bf Theorem~2.1~~}Every binary three numbers  function can be
represented as a superposition of three numbers functions of one
variable.

 A binary three numbers  function is called singular
binary three numbers  function if it's zero in all discrete points
but except one. It's called standard singular three numbers binary
function if non-zero point is in lower-right location. Definitions
for singular three numbers  function of three variables and for
standard singular three numbers function of three variables are
similar.

First we prove that the standard singular three numbers binary
function is a single term one. It's clear the standard singular
binary three numbers function is F above. (0,0,1)
in(0,0,1)$\widetilde{\alpha_{1}}$ in the expression of F is called
raw function. Raw of none-zero point will change if we adjust the
location of `1' in (0,0,1). (-1,-1,0) of(-1,-1,0)
$\widetilde{\alpha_{2}}$ in it is called column fuction.Column of
none-zero point will change if we adjust the location of `0'in
(-1,-1,0). The first (0,0,1) in it is called value function. Value
which may be single-valued or multi-valued or no-valued of non-zero
point will change if we modify `1' in (0,0,1).  Thus we know that
every singular binary three numbers  function can be represented as
a superposition of three numbers functions of one variable.

Because every binary three numbers  function can be transformed to
sum of 9 singular binary three numbers  functions then we get our
theorem.

 So every binary three numbers  function can be represented as:
\begin {align}\label{eq:eps}
\Psi_{2}=\sum_{i=1}^{L}
f_{i}[g_{i1}\widetilde{\alpha_{1}}+g_{i2}\widetilde{\alpha_{2}}]
\end{align}
Here L is not greater than 9. Thus we can express and can construct
a binary three numbers  function by unary three numbers functions.

We can extend all these result to  N numbers function of several
variables.In the decomposition of standard singular binary three
numbers function if we replace raw function (0,0,1) by
(0,0,$\cdots$0,1),column function (-1,-1,0) by (-1,-1,$\cdots$-1,0)
and value function (0,0,1)by(0,0,$\cdots$0,1) ,respectively.Then we
can extend this expression to N numbers functions of two variables.
Situation for N numbers functions of M variables is similar. So we
get conclusion below.

If N$\geq$M+1 a general N numbers  function of M variables can be
decomposed as:
\begin {align}
\psi=\sum_{i=1}^{L}fi\sum _{j=1}^{M}g_{ij}\widetilde{\alpha_{i}}
\end{align}
If N$<$M+1, the number of unary function in expression of singular
discrete  function will be bigger than M+1. For example a standard
singular three numbers function of three variables ¦×~3 can be
represented as:
\begin {align}\label{eq:eps}
\Psi_{3}=(0,0,1)\{(0,0,1)[(0,0,1)\widetilde{\alpha_{1}}+(-1,-1,0)\widetilde{\alpha_{2}}]+(-1,-1,0)\widetilde{\alpha_{3}}\}
\end{align}
Here are more location functions (-1,-1,0) and (0,0,1) than one of
the standard singular binary three numbers  function. Expressions
for singular three numbers function of three variables and for
general three numbers function of three variables are similar to
ones of binary three numbers functions.

All conclusions here are not suit to two numbers function. So we
have:

{\bf Theorem~2.2~~}Every N numbers (N$\geq$3) function of M
variables can be represented as a superposition of N numbers
functions of one variable.

Note we not only prove the existence of representation by
superposition of functions of one variable and give a constructive
procedure. We just only gave the method to decomposing a function
but expression is not the shortest one. Decomposition with terms
being equal to its discrete points is called a trivial
decomposition. Actual terms are more less. Decomposition with less
terms than trivial decomposition is called non-trivial
decomposition. It's an important topic to study non-trivial
decompositions and will not be stated here.

\section{Equations constructed by three numbers functions}

There are $3^{9}$  single-valued binary three numbers  function and
$8^{9}$ ones if they contain many-valued or no-valued ones. How many
equations can we construct with these functions? So many! How many
things need to study about group of the order 3 ? Too poor!  So we
know there are ample mineral resources in this task.

{\bf Theorem~3.1~~}Every algebraic equation constructed by three numbers
functions of two variables can be represented as a superposition of
three numbers functions of one variable.

It's simple to improve it. Solution of any equation is always
function of several variables. By substituting -1,0,1 to the
equation respectively we can get this function easily because field
of definition of it is only three numbers -1,0,1. We can get the
solution expressed by function of one variable by decomposing this
function. That is wonderful that we can construct equations and
solve them freely in a mathematics system! In this paper we just
only solve the equation though there are multitudinous equations:

\qquad\qquad\qquad\qquad\qquad$(x\psi_{1}a)\psi_{3}(x\psi_{2}b)=c$

   Here and in this paper we do'nt write functions of two variables in prefix form  like
$\psi_{3}[\psi_{1}(x,a),\psi_{2}(x,b)]]=c$ for clearness.This equation is called two branches equation. $\psi_{i}$ is parameterized function and can be any one of $8^{9}$ three numbers functions of two
variables. So actually we solve not one equation but a kind of
equation and the method possesses universality.

Assume function $\psi_{1}$,$\psi_{2}$ and $\psi_{3}$ in the two
branches equation is $\Omega_{1}$,$\Omega_{2}$ and $\Omega_{3}$
,respectively:

\qquad \begin{tabular}{|c|c|c|c|} \hline
& -1& 0 & 1\\
\hline
-1&-1& 1& 0  \\
\hline
0&0 & -1 & 1  \\
\hline
1&1 & 0 & -1 \\
\hline
\end{tabular}
{          } \qquad \qquad\begin{tabular}{|c|c|c|c|} \hline
& -1& 0 & 1\\
\hline
-1&0& -1& 1  \\
\hline
0&-1 &0& -1  \\
\hline
1&1 & 1 & 0 \\
\hline
\end{tabular}
{          }
 \qquad \qquad\begin{tabular}{|c|c|c|c|} \hline
& -1& 0 & 1\\
\hline
-1&1& -1& 0  \\
\hline
0&0 &1 & -1  \\
\hline
1&-1 & 0 & 1 \\
\hline
\end{tabular}

When a=b=c=-1 we get numerical equation $[x\Omega_{1}(-1)]\Omega_{3}[x\Omega_{2}(-1)]=-1$. We know only -1 is the solution of this equation by substituting
-1,0,1 to it. So we can know that W(-1,-1,-1)=-1.  By the same way
we can get other values of W(a,b,c). W(a,b,c) can be expressed by
table below.

\begin{tabular}{|c|c|c|c|c|c|c|c|c|c|c|c|}
\hline
c=-1 & -1& 0 & 1 &c=0 & -1& 0 & 1&c=1 & -1& 0 & 1 \\
\hline
-1 &  -1 &N&N&  & 0&N& -1*0*1 & & 1 &  -1*0*1 &N\\
\hline
0 & 1 &  -1*0*1  &N&   & -1 &N&N&  & 0 &N& -1*0*1 \\
\hline
1 & 0 &N&  -1*0*1   &   & 1 &   -1*0*1  &N&  &-1&N&N\\
\hline
\end{tabular}

   In this table the first column indicates the first function number a and
the first row indicates the second function number b and c
indicates the third function number. Decomposing this function of
three variables we get the solution expressed by a superposition of
functions of one variable.

$x=(0,0,-1)\Bigg\{(0,0,1)\Big[(1,0,0)a+(0,-1,-1)b\Big]+(0,-1,-1)c\Bigg\}$

+$(0,0,1)\Bigg\{(0,0,1)\Big[(1,0,0)a+(0,-1,-1)b\Big]+(-1,-1,0)c\Bigg\}$

+$(0,0,N)\Bigg\{(0,0,1)\Big[(1,0,0)a+(-1,0,-1)b\Big]+(0,-1,-1)c\Bigg\}$

+$(0,0,N)\Bigg\{(0,0,1)\Big[(1,0,0)a+(-1,0,-1)b\Big]+(-1,0,-1)c\Bigg\}$

+$(0,0,$-1*0*1$)\Bigg\{(0,0,1)\Big[(1,0,0)a+(-1,0,-1)b\Big]+(-1,-1,0)c\Bigg\}$

+$(0,0,N)\Bigg\{(0,0,1)\Big[(1,0,0)a+(-1,-1,0)b\Big]+(0,-1,-1)c\Bigg\}$

+$(0,0,$-1*0*1$)\Bigg\{(0,0,1)\Big[(1,0,0)a+(-1,-1,0)b\Big]+(-1,0,-1)c\Bigg\}$

+$(0,0,N)\Bigg\{(0,0,1)\Big[(1,0,0)a+(-1,-1,0)b\Big]+(-1,-1,0)c\Bigg\}$

+$(0,0,1)\Bigg\{(0,0,1)\Big[(0,1,0)a+(0,-1,-1)b\Big]+(0,-1,-1)c\Bigg\}$

+$(0,0,-1)\Big]\Bigg\{(0,0,1)\Big[(0,1,0)a+(0,-1,-1)b\Big]+(-1,0,-1)c\Bigg\}$

+$(0,0,$-1*0*1$)\Bigg\{(0,0,1)\Big[(0,1,0)a+(-1,0,-1)b\Big]+(0,-1,-1)c\Bigg\}$

+$(0,0,N)\Bigg\{(0,0,1)\Big[(0,1,0)a+(-1,0,-1)b\Big]+(-1,0,-1)c\Bigg\}$

 +$(0,0,N)\Big]\Bigg\{(0,0,1)\Big[(0,1,0)a+(-1,0,-1)b\Big]+(-1,-1,0)c\Bigg\}$

+$(0,0,N)\Bigg\{(0,0,1)\Big[(0,1,0)a+(-1,-1,0)b\Big]+(0,-1,-1)c\Bigg\}$

+$(0,0,N)\Bigg\{(0,0,1)\Big[(0,1,0)a+(-1,-1,0)b\Big]+(-1,0,-1)c\Bigg\}$

+$(0,0,$-1*0*1$)\Bigg\{(0,0,1)\Big[(0,1,0)a+(-1,-1,0)b\Big]+(-1,-1,0)c\Bigg\}$

+$(0,0,1)\Bigg\{(0,0,1)\Big[(0,0,1)a+(0,-1,-1)b\Big]+(-1,0,-1)c\Bigg\}$

+$(0,0,-1)\Bigg\{(0,0,1)\Big[(0,0,1)a+(0,-1,-1)b\Big]+(-1,-1,0)c\Bigg\}$

+$(0,0,N)\Bigg\{(0,0,1)\Big[(0,0,1)a+(-1,0,-1)b\Big]+(0,-1,-1)c\Bigg\}$

+$(0,0,$-1*0*1$)\Bigg\{(0,0,1)\Big[(0,0,1)a+(-1,0,-1)b\Big]+(-1,0,-1)c\Bigg\}$

+$(0,0,N)\Bigg\{(0,0,1)\Big[(0,0,1)a+(-1,0,-1)b\Big]+(-1,-1,0)c\Bigg\}$

+$(0,0,$-1*0*1$)\Bigg\{(0,0,1)\Big[(0,0,1)a+(-1,-1,0)b\Big]+(0,-1,-1)c\Bigg\}$

+$(0,0,N)\Bigg\{(0,0,1)\Big[(0,0,1)a+(-1,-1,0)b\Big]+(-1,0,-1)c\Bigg\}$

+$(0,0,N)\Big]\Bigg\{(0,0,1)\Big[(0,0,1)a+(-1,-1,0)b\Big]+(-1,-1,0)c\Bigg\}$

 There are 24 but not 27 terms because there are three discrete point with 0 value.

\section{Four special operators}

How to get a new function from a known one? Operator is correspondence between
functions. To give a correspondence between known functions and new functions is to
give an operator. Four special operators mentioned here are easy to
be understood intuitively and are important to solve equations with
parameterized functions however so we must pay attention to them.

{\bf Definition~4.1}Commutation operators. Assume there is an
function of two variables $\psi$, $a_{1}\psi a_{2}=a_{0}$,its
commutation functions
~$\psi$(1,2,0),$\psi$(1,0,2),$\psi$(0,2,1),$\psi$(
\\
2,1,0),$\psi$(2,0,1),$\psi$(0,1,2) will be defined by following
formulas and we introduce commutation operators of one variable C[1,2,0], C[1,0,2],
C[0,2,1], C[2,1,0], C[2,0,1], C[0,1,2] then new functions can be
expressed by $\psi$ and commutation operators.
\begin{subequations}
\begin {align}\label{eq:eps}
a_{1}\psi[1,2,0]a_{2}=a_{0}  \qquad \qquad \psi(1,2,0)=C(1,2,0)(\psi)
\end{align}
\begin {align}
a_{1}\psi[1,0,2]a_{0}=a_{2}\qquad\qquad\psi(1,0,2)=C(1,0,2)(\psi)
\end{align}
\begin {align}
a_{0}\psi[0,2,1]a_{2}=a_{1}\qquad \qquad\psi(0,2,1)=C(0,2,1)(\psi)
\end{align}
\begin {align}
a_{2}\psi[2,1,0] a_{1}=a_{0}\qquad \qquad\psi(2,1,0)=C(2,1,0)(\psi)
\end{align}
\begin {align}
a_{2}\psi[2,0,1] a_{0}=a_{1}\qquad \qquad\psi(2,0,1)=C(2,0,1)(\psi)
\end{align}
\begin {align}
a_{0}\psi[0,1,2]a_{1}=a_{2}\qquad \qquad\psi(0,1,2)=C(0,1,2)(\psi)
\end{align}
\end{subequations}

   Note $\psi$(1,2,0)is $\psi$ itself.

Numbers in brackets indicates new locations of function numbers and
of function result after commutating. That is say original
function doesn't satisfy the new relation gotten by commuting
location of function numbers and of function result but new one
satisfies it. New relation with new function and new location is
equivalent to original one in despite of their forms are different.
For example: if $\Omega$  is the first table below then
~$C(1,0,2)(\Omega)$, $C(0,2,1)(\Omega)$ , $C(2,1,0)(\Omega)$,
$C(2,0,1)(\Omega)$, $C(0,1,2)(\Omega)$will be other tables
,respectively.

\begin{tabular}{|c|c|c|c|} \hline
     & -1& 0 & 1\\
  \hline
   -1&-1& 1& 0  \\
  \hline
  0&0 & -1 & 1  \\
  \hline
  1&1 & 0 & -1 \\
  \hline
\end{tabular}
{          } \qquad \qquad\begin{tabular}{|c|c|c|c|} \hline
& -1& 0 & 1\\
\hline
-1&-1& 1& 0 \\
\hline
0&0 &-1& 1  \\
\hline
1&1 & 0 & -1 \\
\hline
\end{tabular}
{          }
 \qquad \qquad\begin{tabular}{|c|c|c|c|} \hline
 & -1& 0 & 1\\
 \hline
 -1&-1& 0& 1  \\
 \hline
 0&0 &1 & -1  \\
 \hline
 1&1 & -1 & 0 \\
 \hline
\end{tabular}

\begin{tabular}{|c|c|c|c|} \hline
& -1& 0 & 1\\
\hline
-1&-1& 0& 1  \\
\hline
0&1 & -1 & 0  \\
\hline
1&0 & 1 & -1 \\
\hline
\end{tabular}
{          }
  \qquad \qquad\begin{tabular}{|c|c|c|c|} \hline
     & -1& 0 & 1\\
  \hline
   -1&-1& 0& 1  \\
  \hline
  0&0 &1& -1  \\
  \hline
  1&1 & -1 & 0 \\
  \hline
\end{tabular}
{          }
 \qquad \qquad\begin{tabular}{|c|c|c|c|} \hline
     & -1& 0 & 1\\
  \hline
   -1&-1& 0& 1  \\
  \hline
  0&1&-1& 0  \\
  \hline
  1&0 & 1 & -1 \\
  \hline
\end{tabular}

 We can get any combination of function numbers and of
function result for $C(1,0,2)(\Omega)$ by commuting the second
function number and function result for $C(1,2,0)(\Omega)$.
Situations for other commutation functions are similar to it. We
don't limit function at all when we do commutation operator. May be
we get a many-valued function by a not monotonic function or get
an function with no values in some discrete points by a not
surjective function. The same situation may be exists in other
three special operators. We have shown our opinion above. An
mathematics system is extensive and open if it involves solving
equation so it's impossible to limit functions in it. I have ever
tried to limit function in ones of single-valued but failed because
function of many-valued or of no-valued can be introduced from
function of single-valued by special operators. This problem had
troubled me for a long time until I read materials about extension
of group. I known functions of many-valued or of no-valued are not
difficult to be accepted by mathematicians. Commutation operator for
binary functions can be extended to function of many variables.
Showing all commutation functions is integrity in logical and not
all of them will be used in solving equations. There are only two
commutation functions for a unary function:
\begin{subequations}
\begin {equation}
\beta_{e}(a)=a_{0}  \qquad \qquad\qquad   \beta_{e}=\beta
\end{equation}
\begin {equation}
\beta_{t}(a_{0})=a   \qquad \qquad\qquad  \beta_{t}=C(\beta)
\end{equation}
\end{subequations}

{\bf Definition~4.2}Tension-compression operator.Assume there is a
binary function $\psi$ and an unary function $\beta$, $\beta
(a_{1})\psi a_{2}=a_{0}$, we can introduce a new binary function
$\psi _{1}$ by $\psi$ and $\beta$, $\psi _{1}$ will meet the
relation: $a_{1}\psi_{1} a_{2}=a_{0}$, that is say, $a_{1}\psi_{1}
a_{2}=\beta (a_{1})\psi a_{2}$. Introduce a special operator $T_{1}$
to express the relation between  $\psi _{1}$and $\psi$,$\beta$ .
\begin{subequations}
\begin {equation}
   \psi_{1}=\psi T_{1}\beta
\end{equation}
In the same way if $a_{1}\psi \beta (a_{2})=a_{0}$, we can introduce
a new binary function $\psi _{2}$ by $\psi$ and $\beta$, $\psi
_{2}$ will meet the relation: $a_{1}\psi_{2} a_{2}=a_{0}$, that is
say, $a_{1}\psi_{2} a_{2}=a_{1}\psi \beta (a_{2})$. Introduce a
special operator $T_{2}$ to express the relation between  $\psi
_{2}$and $\psi$,$\beta$.
\begin {equation}
  \psi_{2}=\psi T_{2}\beta
\end{equation}
If $a_{1}\psi a_{2}= \beta (a_{0})$,that is say $\beta
^{-1}[a_{1}\psi a_{2}]=a_{0}$,we can introduce a new binary
function $\psi _{0}$ by $\psi$ and $\beta$, $\psi _{0}$ will meet
the relation: $a_{1}\psi_{0} a_{2}=a_{0}$,that is say $a_{1}\psi_{0}
a_{2}=\beta ^{-1}[a_{1}\psi a_{2}]$, and there is $T_{0}$:
\begin {equation}
  \psi_{0}=\psi T_{0}\beta
\end{equation}
\end{subequations}

 For example, (1,-1,0) is an function of one
variable and written in $\gamma$ then ¦¸~$\Omega$$T_{1}$$\gamma$
and ¦¸~$\Omega$$T_{2}$$\gamma$ and ~$\Omega$$T_{0}$$\gamma$ will be

\begin{tabular}{|c|c|c|c|} \hline
     & -1& 0 & 1\\
  \hline
   -1&1& 0& -1  \\
  \hline
  0&-1 & 1 & 0  \\
  \hline
  1&0 & -1 & 1 \\
  \hline
\end{tabular}
{          }
  \qquad \qquad\begin{tabular}{|c|c|c|c|} \hline
     & -1& 0 & 1\\
  \hline
   -1&0& -1& 1  \\
  \hline
  0&1 & 0 & -1  \\
  \hline
  1&-1& 1 & 0 \\
  \hline
\end{tabular}
{          }
 \qquad \qquad\begin{tabular}{|c|c|c|c|} \hline
     & -1& 0 & 1\\
  \hline
   -1&0&-1&1\\
  \hline
  0&1&0&-1\\
  \hline
  1&-1&1&0\\
  \hline
\end{tabular}
respectively.

It's occasional that $\Omega$$T_{2}$$\gamma$ is equal to
$\Omega$$T_{0}$$\gamma$. Only $T_{0}$ will be used in solving
equation.

For an unary function we have only T and $T_{0}$:
\begin{subequations}
\begin {equation}
\beta_{1}T\beta_{2}=\beta_{1}\beta_{2}
\end{equation}
\begin {equation}
\beta_{1}T_{0}\beta_{2}=\beta_{2}^{-1}\beta_{1}
\end{equation}
\end{subequations}

 Note,$\beta_{1}$$\beta_{2}$ means applying$\beta_{2}$ first and
then applying$\beta_{1}$.That is say
\begin {equation}
\\\beta_{1}\beta_{2}(x)=\beta_{1}\Big[\beta_{2}(x)\Big]
\end{equation}
A discrete point for $\beta_{1}$$\beta_{2}$ will be no-valued if it
for any of $\beta_{1}$ or $\beta_{2}$ is no-valued. $\beta_{1}$ and
$\beta_{2}$ will be each other inverse function
if$\beta_{1}$$\beta_{2}$ =e . There are $8^{3}$ three numbers
functions of one variable in which there is always inverse
function for any three numbers function of one variable.

This rule is right for many-valued functions of two variables
because tension-compression operators for functions of two
variables involves actually only composition of two functions of
one variable.

{\bf Definition~4.3}Superposition operator. Assume there are P
functions of many variables$\psi_{k}$(k=1,p),their superposition
function $\psi$ will be:
\begin {equation}
\psi=\sum_{k=1}^{P}\psi_{k}
\end{equation}

Value of $\psi$ will be the sum of value of $\psi_{k}$(k=1,p). This
is a kind of operator by it we can get a new function by several
known functions with same variables. $\psi$ will be no-valued in a
point if any of $\psi_{k}$  is no-valued in this point.
$\psi_{1}+\psi_{2}$ will be many-valued in a point if $\psi_{1}$  is
single-valued and $\psi_{2}$  is many-valued in this point.

{\bf Definition~4.4}Decomposition operator.
\begin {equation}
\psi_{3}=\sum_{i=1}^{27}f_{i}\Bigg\{g_{i4}\Big[g_{i1}(\widetilde{\alpha_{1}})+g_{i2}(\widetilde{\alpha_{2}})\Big]+g_{i3}(\widetilde{\alpha_{3}})\Bigg\}
\end{equation}
We can express the relations between$f_{i}$ or ~$g_{ij}$
and~$\psi_{3}$with special operators $V_{3}$ and $P_{ij}$  and
actually $g_{ij}$ is not change with $\psi_{3}$.
\begin{subequations}
\begin {equation}
f_{i}=V_{i}(\psi_{3})  \qquad   \qquad  \qquad    \qquad (i=1, 27)
\end{equation}
\begin {equation}
g_{ij}=P_{ij}(\psi_{3})        \qquad    \qquad (i=1, 27,j=1,4)
\end{equation}
\end{subequations}

Otherwise there are more than one decomposition for any function of
3 variables but we select only one of them. Correspondence between
$\psi_{3}$ and$f_{i}$ , ~$g_{ij}$ is clear and easy to be gotten. So
decomposition operator is not occult at all.

Please note commutation operator or tension-compression operator or
decomposition operator or superposition operator will be close
within all three numbers functions if they contain ones being
many-valued and no-valued. This is very important and is the
sufficient reason for existing of many-valued functions and
no-valued functions.

So four special operators are very clear and not perplexed at all.

{\bf Definition~4.5}False function of M+K
variables. We  can change an function of M variables to a false one
of (M+K ) variables by adding $o\widetilde{\alpha_{k}}$ in which o
is a zero function and function $\psi$  will not change with K
variables. \qquad\qquad
\begin {equation}
\psi=\sum_{i=1}^{L}fi\sum
_{j=1}^{M}g_{ij}\widetilde{\alpha_{i}}=\sum_{i=1}^{L}fi\Bigg\{\sum
_{j=1}^{M}g_{ij}\widetilde{\alpha_{i}}+\sum
_{k=M+1}^{M+K}o\widetilde{\alpha_{k}}\Bigg\}
\end{equation}

  We can also get false function of (M+K ) variables from one of M variables by $T_{k}o$ (k=M,M+K). For example:

\qquad\qquad\begin{tabular}{|c|c|c|c|c|c|c|c|c|c|c|c|}
  \hline
   c=-1 & -1& 0 & 1 &c=0 & -1& 0 & 1&c=1 & -1& 0 & 1 \\
  \hline
  -1 &  0&1&-1&  & 0,1&N&  -1*0*1  & & 1 & 0&N\\
  \hline
  0 &   0&1&-1&   & 0,1&N&  -1*0*1  &  &  1 & 0&N\\
  \hline
  1 &   0&1&-1 &   &  0,1&N&  -1*0*1 &  & 1 & 0&N\\
  \hline
\end{tabular}

This is a false function of three variables and value of it will
not change with the first variable. Below table is a false function
of two variables.

\qquad\qquad\qquad\qquad\qquad\qquad
\begin{tabular}{|c|c|c|c|}
  \hline
       &-1& 0& 1 \\
  \hline
       -1&1&1&1\\
  \hline
       0&-1&-1&-1\\
  \hline
       1&0&0&0\\
  \hline
\end{tabular}

False function of many variables will be used in solving equations
with parameterized functions.

\section{Formula solution for equations with parameterized functions}

What's an analytic solution or formula solution for an
equation?Formula solution can only contain known parameters or
constants and known parameterized functions or known numerical functions and four kinds
of special operators and we call them valid symbols and all others
invalid ones. This is the standard to verify a formula solution of
an equation. Commutation operators and tension-compression operators
and superposition operators and decomposition operators are the
sufficient condition but not the necessary condition to give formula
solutions of equations. There may be another equivalence set of
operators that can express formula solutions of equations. We will
solve two branches equation as a example below. At the same time we
will solve an equation with digital functions below then we can
understand the procedure more clearly. We must believe that it's not
complex to solve this equation because we have known already the
solution exists surely and only four special operators will be deal
with to get it. We will take any new function met in procedure of
solving the equation as a normal one and will never be puzzled by
its appearance.

\textbf{Step 1:} Decomposing function $\psi_{3}$ as:
 \qquad \qquad\qquad
$$\psi_{3}=\sum_{i=1}^{9}
f_{i}(g_{i1}\widetilde{\alpha_{1}}+g_{i2}\widetilde{\alpha_{2}}) $$

\qquad \qquad\qquad $$\sum_{i=1}^{9}
f_{i}\Big[g_{i1}(x\psi_{1}a)+g_{i2}(x\psi_{2}b)\Big]=c$$
\\$\Omega_{3}=~(0,0,1)\Big[(1,0,0)\widetilde{\alpha_{1}}+(0,-1,-1)\widetilde{\alpha_{2}}\Big]+(0,0,-1)\Big[(1,0,0)\widetilde{\alpha_{1}}+(-1,0,-1)\widetilde{\alpha_{2}}\Big]$
\\$+(0,01)\Big[(0,1,0)\widetilde{\alpha_{1}}+(-1,0,-1)\widetilde{\alpha_{2}}\Big]+(0,0,-1)\Big[(0,1,0)\widetilde{\alpha_{1}}+(-1,-1,0)\widetilde{\alpha_{2}}\Big]$
\\$+(0,0,-1)\Big[(0,0,1)\widetilde{\alpha_{1}}+(0,-1,-1)\widetilde{\alpha_{2}}\Big]+(0,0,1)\Big[(0,0,1)\widetilde{\alpha_{1}}+(-1,-1,0)\widetilde{\alpha_{2}}\Big]$
\\
\\
$(x\Omega_{1}a) \Omega_{3} (x\Omega_{2}b)=$
\\$~(0,0,1)\Big[(1,0,0)(x\Omega_{1}a)+(0,-1,-1)(x\Omega_{2}b)\Big]+(0,0,-1)\Big[(1,0,0)(x\Omega_{1}a)+(-1,0,-1)(
x\Omega_{2}b)\Big]$
\\$+(0,01)\Big[(0,1,0)(x\Omega_{1}a)+(-1,0,-1)( x\Omega_{2}b)\Big]+(0,0,-1)\Big[(0,1,0)(x\Omega_{1}a)+(-1,-1,0)(x\Omega_{2}b)\Big]$
\\$+(0,0,-1)\Big[(0,0,1)(x\Omega_{1}a)+(0,-1,-1)(x\Omega_{2}b)\Big]+(0,0,1)\Big[(0,0,1)(x\Omega_{1}a)+(-1,-1,0)(x\Omega_{2}b)\Big]$
=c

\textbf{Step 2:} By tension-compression of$g_{i1},g_{i2}$ we have:

\qquad\qquad$$\sum_{i=1}^{9}f_{i}\Big[x(\psi_{1}T_{0}g_{i1}^{-1})a+x(\psi_{2}T_{0}g_{i2}^{-1})b\Big]=c$$

Note,$(\psi_{1}T_{0}g_{i1}^{-1})$ in $x(\psi_{1}T_{0}g_{i1}^{-1})a$
and $(\psi_{2}T_{0}g_{i2}^{-1})$ in $x(\psi_{2}T_{0}g_{i2}^{-1})b$
are two functions of two variables.

$(x\Omega_{1}a) \Omega_{3} (x\Omega_{2}b)=~$
$(0,0,1)\Bigg\{x\Big[\Omega_{1}T_{0}(1,0,0)^{-1}\Big]a+x\Big[\Omega_{2}T_{0}(0,-1,-1)^{-1}\Big]b\Bigg\}+$

$(0,0,-1)\Bigg\{x\Big[\Omega_{1}T_{0}(1,0,0)^{-1}\Big]a+x\Big[\Omega_{2}T_{0}(-1,0,-1)^{-1}\Big]b\Bigg\}+$
$(0,01)\Bigg\{x\Big[\Omega_{1}T_{0}(0,1,0)^{-1}\Big]a$

$+x\Big[\Omega_{2}T_{0}(-1,0,-1)^{-1}\Big]b\Bigg\}+$
$(0,0,-1)\Bigg\{x\Big[\Omega_{1}T_{0}(0,1,0)^{-1}\Big]a+x\Big[\Omega_{2}T_{0}(-1,-1,0)^{-1}\Big]b\Bigg\}+$

$(0,0,-1)\Bigg\{x\Big[\Omega_{1}T_{0}(0,0,1)^{-1}\Big]a+x\Big[\Omega_{2}T_{0}(0,-1,-1)^{-1}\Big]b\Bigg\}+$
$(0,0,1)\Bigg\{x\Big[\Omega_{1}T_{0}(0,0,1)^{-1}\Big]a+$

$x\Big[\Omega_{2}T_{0}(-1,-1,0)^{-1}\Big]b\Bigg\}=c$

$\Omega_{1}T_{0}(1,0,0)^{-1}$, $\Omega_{1}T_{0}(0,1,0)^{-1}$,
$\Omega_{1}T_{0}(0,0,1)^{-1}$ is

\begin{tabular}{|c|c|c|c|}
\hline
     & -1& 0 & 1\\
  \hline
   -1&1& 0& 0  \\
  \hline
  0&0 &1 & 0 \\
  \hline
  1&0 & 0 & 1 \\
  \hline
\end{tabular}
{}
  \qquad \qquad\qquad \begin{tabular}{|c|c|c|c|} \hline
     & -1& 0 & 1\\
  \hline
   -1&0& 0& 1  \\
  \hline
  0&1 &0 & 0 \\
  \hline
  1&0 & 1 & 0 \\
  \hline
\end{tabular}
{}
 \qquad \qquad\begin{tabular}{|c|c|c|c|} \hline
     & -1& 0 & 1\\
  \hline
   -1&0& 1& 0  \\
  \hline
  0&0 &0 & 1 \\
  \hline
  1&1 & 0 & 0 \\
  \hline
\end{tabular}

$\Omega_{2}T_{0}(0,-1,-1)^{-1}$, $\Omega_{2}T_{0}(-1,0,-1)^{-1}$,
$\Omega_{2}T_{0}(-1,-1,0)^{-1}$ is

\begin{tabular}{|c|c|c|c|} \hline
     & -1& 0 & 1\\
  \hline
   -1&-1& 0& -1  \\
  \hline
  0&0& -1& 0 \\
  \hline
  1&-1& -1& -1\\
  \hline
\end{tabular}
{}
  \qquad \qquad\begin{tabular}{|c|c|c|c|} \hline
     & -1& 0 & 1\\
  \hline
   -1&0& -1& -1  \\
  \hline
  0&-1& 0& -1 \\
  \hline
  1&-1& -1& 0\\
  \hline
\end{tabular}
{}
 \qquad \qquad\begin{tabular}{|c|c|c|c|} \hline
     & -1& 0 & 1\\
  \hline
   -1&-1& -1& 0  \\
  \hline
  0&-1& -1& -1 \\
  \hline
  1&0& 0& -1\\
  \hline
\end{tabular}
respectively.

\textbf{Step 3:} Changing $\psi_{1}T_{0}g_{i1}^{-1}$ by $T_{3}o$ to
get a false function of three variables
$\psi_{1}T_{3}oT_{0}g_{i1}^{-1}$ in which variable c is a false one
and Changing $\psi_{2}T_{0}g_{i2}^{-1}$ by $T_{2}$o to get a false
function of three variables $\psi_{2}T_{2}oT_{0}g_{i2}^{-1}$ in
which variable b is a false one ,respectively. Adding them to get a
real function of three variables $\psi_{i3}$. This is the
application of tension-compression operator in solving equation.

\qquad\qquad$\psi_{i3}=\psi_{1}T_{3}oT_{0}g_{i1}^{-1}+\psi_{2}T_{2}oT_{0}g_{i2}^{-1}\qquad(i=1,9)$

$\theta_{1}=\Omega_{1}T_{3}oT_{0}(1,0,0)^{-1}+\Omega_{2}T_{2}oT_{0}(0,-1,-1)^{-1}$
is:

\begin{tabular}{|c|c|c|c|c|c|c|c|c|c|c|c|}
  \hline
   c=-1 & -1& 0 & 1 &c=0 & -1& 0 & 1&c=1 & -1& 0 & 1 \\
  \hline
  -1 & 0 & -1 & -1 &   & 1 & 0 & 0 &  & 0 & -1 & -1  \\
  \hline
  0 & 0 & 1 & 0&   & -1 & 0 & -1&  & 0 & 1 & 0 \\
  \hline
  1 & -1 & -1 & 0 &   & -1 & -1 & 0 &  & -1 & -1 & 0  \\
  \hline
\end{tabular}

$\theta_{2}=\Omega_{1}T_{3}oT_{0}(1,0,0)^{-1}+\Omega_{2}T_{2}oT_{0}(-1,0,-1)^{-1}$
is:

\begin{tabular}{|c|c|c|c|c|c|c|c|c|c|c|c|}
  \hline
   c=-1 & -1& 0 & 1 &c=0 & -1& 0 & 1&c=1 & -1& 0 & 1 \\
  \hline
  -1 & 1& 0 & 0&   & 0 & -1 & -1 &  & 0 & -1 & -1  \\
  \hline
  0 & -1 & 0 & -1&   & 0 & 1 & 0&  & -1 & 0& -1 \\
  \hline
  1 & -1 & -1 & 0 &   & -1 & -1 & 0 &  & 0 & 0& 1  \\
  \hline
\end{tabular}

$\theta_{3}=\Omega_{1}T_{3}oT_{0}(0,1,0)^{-1}+\Omega_{2}T_{2}oT_{0}(-1,0,-1)^{-1}$
is:

\begin{tabular}{|c|c|c|c|c|c|c|c|c|c|c|c|}
  \hline
   c=-1 & -1& 0 & 1 &c=0 & -1& 0 & 1&c=1 & -1& 0 & 1 \\
  \hline
  -1 & 0 &0 & 1 &   & -1 & -1 & 0 &  & -1 & -1 & 0  \\
  \hline
  0 & 0 & -1 & -1&   & 1 & 0 & 0&  & 0 & -1 & -1 \\
  \hline
  1 & -1 & 0 & -1 &  & -1 & 0 & -1 &  & 0 & 1 & 0  \\
  \hline
\end{tabular}

$\theta_{4}=\Omega_{1}T_{3}oT_{0}(0,1,0)^{-1}+\Omega_{2}T_{2}oT_{0}(-1,-1,0)^{-1}$
is:

\begin{tabular}{|c|c|c|c|c|c|c|c|c|c|c|c|}
  \hline
   c=-1 & -1& 0 & 1 &c=0 & -1& 0 & 1&c=1 & -1& 0 & 1 \\
  \hline
  -1 &-1  & -1 & 0&  & -1 & -1 & 0 &  & 0 & 0 & 1  \\
  \hline
  0 & 0 & -1  & -1 &   & 0 & -1  & -1&  & 0 & -1  & -1  \\
  \hline
  1 & 0 & 1 & 0 &   & 0 & 1 & 0 &  & -1 & 0 & -1   \\
  \hline
\end{tabular}

$\theta_{5}=\Omega_{1}T_{3}oT_{0}(0,0,1)^{-1}+\Omega_{2}T_{2}oT_{0}(0,-1,-1)^{-1}$
is:

\begin{tabular}{|c|c|c|c|c|c|c|c|c|c|c|c|}
  \hline
   c=-1 & -1& 0 & 1 &c=0 & -1& 0 & 1&c=1 & -1& 0 & 1 \\
  \hline
  -1 & -1 & 0 & -1 &   & 0 & 1 & 0 & & -1& 0 & -1  \\
  \hline
  0 & 0 & 0 & 1&  & -1 & -1 & 0&  & 0 & 0 & 1 \\
  \hline
  1 & 0 & -1 & -1 &   & 0 & -1 & -1 &  & 0 & -1 & -1  \\
  \hline
\end{tabular}

$\theta_{6}=\Omega_{1}T_{3}oT_{0}(0,0,1)^{-1}+\Omega_{2}T_{2}oT_{0}(-1,-1,0)^{-1}$
is:

\begin{tabular}{|c|c|c|c|c|c|c|c|c|c|c|c|}
  \hline
   c=-1 & -1& 0 & 1 &c=0 & -1& 0 & 1&c=1 & -1& 0 & 1 \\
  \hline
  -1 & -1 & 0 & -1 &   & -1 & 0 & -1 &  & 0 & 1 & 0  \\
  \hline
  0 & -1& -1 & 0&   & -1 & -1 & 0& & -1 & -1 & 0 \\
  \hline
  1 & 1 & 0 & 0 &  & 1 & 0& 0 & & 0 & -1 & -1  \\
  \hline
\end{tabular}
\\

\textbf{Step 4: }Changing $\psi_{i3}$ by $T_{0}f_{i}^{-1}$we get:

\qquad\qquad \qquad \qquad \qquad\qquad \qquad
$\psi_{i4}=\psi_{i3}T_{0}f_{i}^{-1}\qquad(i=1,9)$

$\theta_{1}T_{0}(0,0,1)^{-1}$ is:

\begin{tabular}{|c|c|c|c|c|c|c|c|c|c|c|c|}
  \hline
   c=-1 & -1& 0 & 1 &c=0 & -1& 0 & 1&c=1 & -1& 0 & 1 \\
  \hline
  -1 & 0 & 0 & 0 &   & 1 & 0 & 0 &  & 0 & 0 & 0  \\
  \hline
  0 & 0 & 1 & 0&   & 0 & 0 & 0&  & 0 & 1 & 0 \\
  \hline
  1 & 0 & 0 & 0 &   & 0 & 0& 0 &  & 0 & 0 & 0  \\
  \hline
\end{tabular}

$\theta_{2}T_{0}(0,0,-1)^{-1}$ is:

\begin{tabular}{|c|c|c|c|c|c|c|c|c|c|c|c|}
  \hline
   c=-1 & -1& 0 & 1 &c=0 & -1& 0 & 1&c=1 & -1& 0 & 1 \\
  \hline
  -1 & -1 & 0 & 0&   & 0 & 0& 0 &  & 0 & 0 & 0  \\
  \hline
  0 & 0 & 0 & 0&  & 0 & -1  & 0&  & 0 & 0& 0 \\
  \hline
  1 & 0 & 0 & 0 &   & 0 & 0 & 0 &  & 0 & 0& -1  \\
  \hline
\end{tabular}

$\theta_{3}T_{0}(0,0,1)^{-1}$ is:

\begin{tabular}{|c|c|c|c|c|c|c|c|c|c|c|c|}
  \hline
   c=-1 & -1& 0 & 1 &c=0 & -1& 0 & 1&c=1 & -1& 0 & 1 \\
  \hline
  -1 & 0 &0 & 1 &  & 0 & 0 & 0 & & 0 & 0& 0  \\
  \hline
  0 & 0 & 0 & 0&   & 1 & 0 & 0& & 0 & 0 & 0 \\
  \hline
  1 & 0 & 0 & 0 &   & 0 & 0 & 0 &  & 0 & 1 & 0  \\
  \hline
\end{tabular}

$\theta_{4}T_{0}(0,0,-1)^{-1}$ is:

\begin{tabular}{|c|c|c|c|c|c|c|c|c|c|c|c|}
  \hline
   c=-1 & -1& 0 & 1 &c=0 & -1& 0 & 1&c=1 & -1& 0 & 1 \\
  \hline
  -1 &0  & 0 & 0&   & 0 & 0 & 0 &  & 0 & 0 & -1  \\
  \hline
  0 & 0 & 0  & 0 &   & 0 & 0  & 0&  & 0 & 0  & 0  \\
  \hline
  1 & 0 & -1& 0 &  & 0 & -1 & 0 &  & 0 & 0 & 0   \\
  \hline
\end{tabular}

$\theta_{5}T_{0}(0,0,-1)^{-1}$ is:

\begin{tabular}{|c|c|c|c|c|c|c|c|c|c|c|c|}
  \hline
   c=-1 & -1& 0 & 1 &c=0 & -1& 0 & 1&c=1 & -1& 0 & 1 \\
  \hline
  -1 & 0 & 0 & 0 &   & 0 & -1 & 0 &  & 0& 0 & 0 \\
  \hline
  0 & 0 & 0 & -1&  & 0 & 0 & 0&  & 0 & 0 & -1 \\
  \hline
  1 & 0 & 0 & 0 &   & 0 & 0 & 0 &  & 0 & 0 & 0  \\
  \hline
\end{tabular}

$\theta_{6}T_{0}(0,0,1)^{-1}$ is:

\begin{tabular}{|c|c|c|c|c|c|c|c|c|c|c|c|}
  \hline
   c=-1 & -1& 0 & 1 &c=0 & -1& 0 & 1&c=1 & -1& 0 & 1 \\
  \hline
  -1 & 0 & 0 & 0 &   & 0 & 0 & 0 &  & 0 & 1 & 0  \\
  \hline
  0 & 0& 0 & 0&   & 0 & 0 & 0&  & 0 & 0 & 0 \\
  \hline
  1 & 1 & 0 & 0 &   & 1 & 0& 0 &  & 0 & 0 & 0 \\
  \hline
\end{tabular}

\textbf{Step 5:} To sum $\psi_{i4}$  we get:

\qquad  \qquad  \qquad \qquad \qquad  \qquad
$$\psi_{5}=\sum_{i=1}^{9} \psi_{i4} $$

Original equation will be:

 \qquad  \qquad  \qquad \qquad
\qquad  \qquad $\psi_{5}(x,a,b)=c$

$\theta_{7}=\theta_{1}T_{0}(0,0,1)^{-1}+\theta_{2}T_{0}(0,0,-1)^{-1}+\theta_{3}T_{0}(0,0,1)^{-1}+\theta_{4}T_{0}(0,0,-1)^{-1}$

+$\theta_{5}T_{0}(0,0,-1)^{-1}+\theta_{6}T_{0}(0,0,1)^{-1}$ is:

\begin{tabular}{|c|c|c|c|c|c|c|c|c|c|c|c|}
\hline
c=-1 & -1& 0 & 1 &c=0 & -1& 0 & 1&c=1 & -1& 0 & 1 \\
\hline
-1 &  -1 & 0 & 1 &   & 1 &  -1 & 0 & & 0 & 1 &  -1  \\
\hline
0 & 0& 1 &  -1&  & 1 &  -1 & 0&  & 0 & 1 &  -1 \\
\hline
1 & 1 &  -1 & 0 &   & 1 &  -1& 0 &  & 0 & 1 &  -1 \\
\hline
\end{tabular}

\textbf{Step6: }By commutation operator we get:

\qquad \qquad $x=\Big[C(2,3,0,1) \psi_{5}\Big](a,b,c)=W(a,b,c)$
\\
$C(2,3,0,1)\theta_{7}$is:

\begin{tabular}{|c|c|c|c|c|c|c|c|c|c|c|c|}
\hline
c=-1 & -1& 0 & 1 &c=0 & -1& 0 & 1&c=1 & -1& 0 & 1 \\
\hline
-1 &  -1 &N&N&  & 0&N& -1*0*1 & & 1 & -1*0*1&N\\
\hline
0 & 1 & -1*0*1 &N&   & -1 &N&N&  & 0 &N&-1*0*1\\
\hline
1 & 0 &N& -1*0*1  &   & 1 &  -1*0*1 &N&  &-1&N&N\\
\hline
\end{tabular}

It's not oddball there are many-valued discrete points or no-valued
discrete points. Not all commutation operators are used in solving
equation.

\textbf{Step 7:} by decomposition operator we get:

 $$x=\sum_{k=1}^{27}u_{k}\Bigg\{v_{k4}\Big[v_{k1}(a)+v_{k2}(b)\Big]+v_{k3}(c)\Bigg\}$$
$$=\sum_{k=1}^{27}(V_{k}W)\Bigg\{(P_{k4}W)\Big[(P_{k1}W)(a)+(P_{k2}W)(b)\Big]+(P_{k3}W)(c)\Bigg\}$$

we replace logogram symbols by complete ones.

$$x=\sum_{k=1}^{27}V_{k}\Big[C(2,3,0,1)\Big(\sum_{i=1}^{9}\big\{\big[\psi_{1}T_{3}oT_{0}(P_{i1}\psi_{3})^{-1}+\psi_{2}T_{2}oT_{0}(P_{i2}\psi_{3})^{-1}\big]T_{0}(V_{i}\psi_{3})^{-1}\big\}\Big)\Big]$$

$$\Bigg[P_{k4}\Big[C(2,3,0,1)\Big(\sum_{i=1}^{9}\big\{\big[\psi_{1}T_{3}oT_{0}(P_{i1}\psi_{3})^{-1}+\psi_{2}T_{2}oT_{0}(P_{i2}\psi_{3})^{-1}\big]T_{0}(V_{i}\psi_{3})^{-1}\big\}\Big)\Big]$$

$$\Bigg(P_{k1}\Big[C(2,3,0,1)\Big(\sum_{i=1}^{9}\big\{\big[\psi_{1}T_{3}oT_{0}(P_{i1}\psi_{3})^{-1}+\psi_{2}T_{2}oT_{0}(P_{i2}\psi_{3})^{-1}\big]T_{0}(V_{i}\psi_{3})^{-1}\big\}\Big)\Big](a)$$

$$+P_{k2}\Big[C(2,3,0,1)\Big(\sum_{i=1}^{9}\big\{\big[\psi_{1}T_{3}oT_{0}(P_{i1}\psi_{3})^{-1}+\psi_{2}T_{2}oT_{0}(P_{i2}\psi_{3})^{-1}\big]T_{0}(V_{i}\psi_{3})^{-1}\big\}\Big)\Big](b)\Bigg)$$

$$+P_{k3}\Big[C(2,3,0,1)\Big(\sum_{i=1}^{9}\big\{\big[\psi_{1}T_{3}oT_{0}(P_{i1}\psi_{3})^{-1}+\psi_{2}T_{2}oT_{0}(P_{i2}\psi_{3})^{-1}\big]T_{0}(V_{i}\psi_{3})^{-1}\big\}\Big)\Big](c)\Bigg]$$

Actually location functions $P_{ij}\psi_{k}$ do not change with
$\psi_{k}$ and can be called constant  functions. Solution of
equation with function $\Omega_{1}$,$\Omega_{2}$ and $\Omega_{3}$
has been given already above by getting a function of many variables
. Giving the procedure of it is just only make the method clearer.
We deal with the function of three variables in solving this
equation. Can we avoid to use it in the procedure? Never!  In
history one reason to introduce complex number is that we have to
deal with complex number even if three roots of a cubic equation are
all real number. It's the most important that we have gotten the
solution expressed by function of one variable however.

\section{Composition of special operators}

There are 10 compositions for commutation  operators and
tension-compression operators and superposition operators and
decomposition operators as bellow tables:
\\
table1
\\
\begin{tabular}{|c|c|c|c|c|}
  \hline
     & commutation & tension-compression& superposition & decomposition \\
   \hline
  commutation      & 1 &  & &  \\
  \hline
  tension-compression & 2& 3&   &   \\
   \hline
  superposition & 4 & 5& 6&  \\
  \hline
  decomposition & 7& 8& 9 & 10  \\
  \hline
\end{tabular}

We will mention them below. Here we give only results of binary
function and they can be extended to functions of many variables
easily.
\\
Composition 1 commutation and commutation:

see table 2
\\
Composition 2 tension-compression and commutation:

see table 3
\\
Composition 3 tension-compression and tension-compression:

see table4
\\
Composition 4 superposition and commutation:

is equal to commutation and superposition for commutation C(2,1,0):
\begin {equation}
C(2,1,0)(\sum_{k=1}^{H}\psi_{k})=\sum_{k=1}^{H}\Big[C(2,1,0)(\psi_{k})\Big]
\qquad
\end{equation}
It will be complex for commutation C(0,2,1) and commutation
C(1,0,2).
\\
Composition 5 superposition and tension-compression:

is equal to tension-compression and superposition for
tension-compression $T_{1}$ and $T_{2}$:
\begin{subequations}
\begin {equation}
(\sum_{k=1}^{H}\psi_{k})T_{1}\beta=\sum_{k=1}^{H}(\psi_{k}T_{1}\beta)
\qquad \qquad \qquad
\end{equation}
\begin {equation}
(\sum_{k=1}^{H}\psi_{k})T_{2}\beta=\sum_{k=1}^{H}(\psi_{k}T_{2}\beta)
\qquad \qquad \qquad
\end{equation}
\end{subequations}
is complex for tension-compression $T_{0}$.
\\
Composition 6 superposition and superposition:

It is very simple.
\\
Composition 7 decomposition - commutation:

Value functions will hold the line and location functions will
exchange for commutation C(2,1,0)
\begin{subequations}
\begin {equation}
V_{i}\Big[C(2,1,0)(\psi) \Big]=V_{i}(\psi) \qquad (i=1,L)
\end{equation}
\begin {equation}
P_{i1}\Big[C(2,1,0)(\psi) \Big]=P_{i2}(\psi) \qquad (i=1,L)
\end{equation}
\begin {equation}
P_{i2}\Big[C(2,1,0)(\psi) \Big]=P_{i1}(\psi) \qquad (i=1,L)
\end{equation}
\end{subequations} is complex for commutation C(0,2,1) and
commutation C(1,0,2).
\\
Composition 8 decomposition and tension-compression:

Value functions will hold the line and location functions will be
acted by $T_{1}\beta$ or $T_{2}\beta$ for tension-compression
$T_{1}$ and  $T_{2}$.
\begin{subequations}
\begin {equation}
V_{i}(\psi T_{j}\beta)=V_{i}\psi \qquad (i=1,L\qquad j=1,2)
\end{equation}
\begin {equation}
P_{ij}(\psi T_{j}\beta)=(P_{ij}\psi)T_{j}\beta  \qquad (i=1,L\qquad
j=1,2)
\end{equation}
\end{subequations}

is complex for $T_{0}$. But there is relation between value
functions of $\psi$ and of $\psi$ acted by $T_{0}$ if it's a
trivial decomposition.
\begin {equation}
V_{i}(\psi T_{0}\beta)=(V_{i}\psi)T_{0}\beta \qquad (i=1,L)
\end{equation}
This relation is very important.
\\
Composition 9 decomposition and superposition:

Value functions will be composition of value functions and
location functions will be any location functions.
\begin{subequations}
\begin {equation}
V_{i}(\sum_{k=1}^{H} \psi _{k})=\sum_{k=1}^{H}V_{i}(\psi _{k})
\qquad (i=1,L )
\end{equation}
\begin {equation}
P_{ij}(\sum _{k=1}^{H}\psi _{i})=P_{ij}(\psi _{k}) \qquad (i=1,L
\qquad j=1,2 )
\end{equation}
\end{subequations}
\\
Composition 10 decomposition and decomposition:

None.

Law and composition of special operators can be extended to high
degree operators in form.
\\
Table2 commutation and commutation

\qquad  \qquad\begin{tabular}{|c|c|c|c|c|c|c|}
  \hline
          & C(1,2,0)& C(1,0,2) & C(0,2,1) & C(2,1,0)& C(2,0,1)& C(0,1,2)\\
   \hline
  C(1,2,0)& C(1,2,0) & C(1,0,2)& C(0,2,1) & C(2,1,0) & C(2,0,1)& C(0,1,2) \\
  \hline
  C(1,0,2)& C(1,0,2) & C(1,2,0)& C(2,0,1) & C(0,1,2) & C(0,2,1)& C(2,1,0) \\
   \hline
  C(0,2,1)& C(0,2,1) & C(0,1,2)& C(1,2,0) & C(2,0,1) & C(2,1,0)& C(1,0,2) \\
  \hline
  C(2,1,0)& C(2,1,0) & C(2,0,1)& C(0,1,2) & C(1,2,0) & C(1,0,2)& C(0,2,1) \\
  \hline
  C(2,0,1)& C(2,0,1) & C(2,1,0)& C(1,0,2) & C(0,2,1) & C(0,1,2)& C(1,2,0) \\
  \hline
  C(0,1,2)& C(0,1,2) & C(0,2,1)& C(2,1,0) & C(1,0,2) & C(1,2,0)& C(2,0,1) \\
 \hline
\end{tabular}
\\
Table3 tension-compression and commutation

\begin{tabular}{|c|c|c|c|c|c|c|}
\hline
& C(1,2,0)& C(1,0,2) & C(0,2,1) & C(2,1,0)& C(2,0,1)& C(0,1,2) \\
\hline
$T_{1}\beta $ & $T_{1}\beta $ & $C(1,0,2)T_{1}\beta $& $C(0,2,1)T_{0}\beta $ & $C(2,1,0)T_{2}\beta $ & $C(2,0,1)T_{0}\beta $& $C(0,1,2)T_{2}\beta $ \\
\hline
$T_{2}\beta $ & $T_{2}\beta $ & $C(1,0,2)T_{0}\beta $& $C(0,2,1)T_{2}\beta $ & $C(2,1,0)T_{1}\beta $ & $C(2,0,1)T_{1}\beta $& $C(0,1,2)T_{0}\beta $ \\
\hline
$T_{0}\beta $ & $T_{0}\beta $ & $C(1,0,2)T_{2}\beta $& $C(0,2,1)T_{1}\beta $ & $C(2,1,0)T_{0}\beta $ & $C(2,0,1)T_{2}\beta $& $C(0,1,2)T_{1}\beta $ \\
\hline
\end{tabular}
\\
table4 tension-compression and tension-compression

\qquad\qquad\qquad\begin{tabular}{|c|c|c|c|c|c|c|} \hline
&$T_{1}\beta_{2}$      &$T_{2}\beta_{2}$          & $T_{0}\beta_{2}$      \\
\hline
$T_{1}\beta_{1}$ &$T_{1}(\beta_{1}\beta_{2})$&$(T_{2}\beta_{2})T_{1}\beta_{1}$&($T_{0}\beta_{2})T_{1}\beta_{1}$ \\
\hline
$T_{2}\beta_{1}$ &$(T_{1}\beta_{2})T_{2}\beta_{1}$&$T_{2}(\beta_{1}\beta_{2})$&$(T_{0}\beta_{2})T_{2}\beta_{1}$\\
\hline
$T_{0}\beta_{1}$ &$(T_{1}\beta_{2})T_{0}\beta_{1}$&$(T_{2}\beta_{2})T_{0}\beta_{1}$&$T_{0}(\beta_{1}\beta_{2})$\\
\hline
\end{tabular}

All of them are easy to be validated by readers.

\section{Extend results to discrete operators}

Now we extend results about discrete functions  to discrete
operators. We limit the field of definition and range of operators
within three discrete functions -e=(1,0,-1),o=(0,0,0),e=(-1,0,1) for
simplicity.

{\bf Definition~7.1}Assume there are three numbers functions
-e=(1,0,-1),o=(0,0,0) and e=(-1,0,1)we let A=$\{$-e,0,e$\}$ and
define three functions operator of one variable $S_{1}$ as

\qquad \qquad\qquad \qquad\qquad \qquad $S_{1}$:A$\longrightarrow$A

define three functions operator of two variables $S_{2}$ as

\qquad \qquad\qquad \qquad\qquad \qquad
$S_{2}$:$A^{2}$$\longrightarrow$A£¬

define three functions operator of three variables $S_{3}$ as

\qquad \qquad\qquad \qquad\qquad
\qquad$S_{3}$:$A^{3}$$\longrightarrow$A

There are $3^{3}$ single-valued three functions operators of one
variable and $8^{3}$ ones if they contain many-valued or
no-valued.There are $3^{9}$ single-valued three functions operators
of two variables and $8^{9}$ ones if they contain many-valued or
no-valued.There are $3^{27}$ single-valued three functions operators
of three variable and $8^{27}$ ones if they contain many-valued or
no-valued.

Functions will be partitioned by the symbol '*'for many-valued point
and no-valued point will be indicated by 'N'.

'+'operation will be expressed as:

\qquad \qquad\qquad \qquad\qquad \qquad\qquad
\begin{tabular}{|c|c|c|c|}
  \hline
       &-1& 0& 1 \\
  \hline
       -1&1&-1&0\\
  \hline
       0&-1&0&1\\
  \hline
       1&0&1&-1\\
  \hline
\end{tabular}

'+'operator will be expressed as:

\qquad \qquad\qquad \qquad\qquad\qquad \qquad
\begin{tabular}{|c|c|c|c|}
  \hline
       &-e& o& e \\
  \hline
       -e&e&-e&o\\
  \hline
       o&-e&o&e\\
  \hline
       e&o&e&-e\\
  \hline
\end{tabular}

Compare two tables we know -e,o,e in discrete operators system is
like -1,0,1 in discrete functions system ,respectively. We can also
introduce concepts of singular three functions operator and standard
singular three functions operator.

A standard singular three functions operator of two variables can be
expressed by table:

\qquad \qquad\qquad\qquad\qquad\qquad \qquad
\begin{tabular}{|c|c|c|c|}
  \hline
       &-e& o& e \\
  \hline
       -e&o&o&o\\
  \hline
       o&o&o&o\\
  \hline
       e&o&o&e\\
  \hline
\end{tabular}

It can be represented as a superposition of three functions
operators of one variable:

 \qquad \qquad \qquad\qquad
$G=(o,o,e)\Big[(o,o,e)\widetilde{\beta_{1}}+(-e,-e,o)\widetilde{\beta_{2}}\Big]$

By the same reason for three numbers function we know a standard
singular binary three functions operator can be represented as a
superposition of unary three functions operators and so does a
general singular binary three functions operator. A general binary
three functions operator can be expressed to sum of 9 singular
binary three functions operators so we have

{\bf Theorem~7.1~~}Every binary three functions operator can be
represented as a superposition of three functions operators of one
variable.

A standard singular three functions operator of three variables
$\phi_{3}$ can be represented as:

$\phi_{3}=(o,o,e)\Bigg\{(o,o,e)\Big[(o,o,e)(\widetilde{\beta_{1}})+(-e,-e,o)(\widetilde{\beta_{2}})\Big]+(-e,-e,o)(\widetilde{\beta_{3}})\Bigg\}$

{\bf Theorem~7.2~~}Every three functions operator of two or of three
variables can be represented as a superposition of three functions
operators of one variable.

All conclusions here are not suit to discrete 2 operator.

There are great number of operator equations constructed by $8^{9}$
operators of two variables.

{\bf Theorem~7.3~~}Every operator equation constructed by three functions
operators of two variables can be give formula solution represented
as a superposition of three functions operators of one variable.

Although there are many operator equations we give formula solution
for only double branches operator equation with digital operators
and with parameterized operators.

\qquad \qquad \qquad \qquad      \qquad \qquad
$(y\phi_{1}f)\phi_{3}(y\phi_{2}g)=h$

Assume~$\phi_{1}$,$\phi_{2}$,$\phi_{3}$ is
~$\Theta_{1},\Theta_{2},\Theta_{3}$ as below,respectively:

\begin{tabular}{|c|c|c|c|} \hline
     & -e& o & e\\
  \hline
   -e&-e& e& o  \\
  \hline
  o&o & -e & e  \\
  \hline
  e&e & o & -e \\
  \hline
\end{tabular}
{          }
  \qquad \qquad\begin{tabular}{|c|c|c|c|} \hline
     & -e& o & e\\
  \hline
   -e&o& -e& e  \\
  \hline
  o&-e &o& -e  \\
  \hline
  e&e & e & o \\
  \hline
\end{tabular}
{          }
 \qquad \qquad\begin{tabular}{|c|c|c|c|} \hline
     & -e& o & e\\
  \hline
   -e&e& -e& o  \\
  \hline
  o&o &e & -e  \\
  \hline
  e&-e & o & e \\
  \hline
\end{tabular}

Solution expressed by superposition of operators of one variable
will be:

$y=(o,o,-e)\Bigg\{(o,o,e)\Big[(e,o,o)f+(o,-e,-e)g\Big]+(o,-e,-e)h\Bigg\}$

+$(o,o,e)\Bigg\{(o,o,e)\Big[(e,o,o)f+(o,-e,-e)g\Big]+(-e,-e,o)h\Bigg\}$

+$(o,o,N)\Bigg\{(o,o,e)\Big[(e,o,o)f+(-e,o,-e)g\Big]+(o,-e,-e)h\Bigg\}$

+$(o,o,N)\Bigg\{(o,o,e)\Big[(e,o,o)f+(-e,o,-e)g\Big]+(-e,o,-e)h\Bigg\}$

+$(o,o,$-e*o*e$)\Bigg\{(o,o,e)\Big[(e,o,o)f+(-e,o,-e)g\Big]+(-e,-e,o)h\Bigg\}$

+$(o,o,N)\Bigg\{(o,o,e)\Big[(e,o,o)f+(-e,-e,o)g\Big]+(o,-e,-e)h\Bigg\}$

+$(o,o,$-e*o*e$)\Bigg\{(o,o,e)\Big[(e,o,o)f+(-e,-e,o)g\Big]+(-e,o,-e)h\Bigg\}$

+$(o,o,N)\Bigg\{(o,o,e)\Big[(e,o,o)f+(-e,-e,o)g\Big]+(-e,-e,o)h\Bigg\}$

+$(o,o,e)\Bigg\{(o,o,e)\Big[(o,e,o)f+(o,-e,-e)g\Big]+(o,-e,-e)h\Bigg\}$

+$(o,o,-e)\Big]\Bigg\{(o,o,e)\Big[(o,e,o)f+(o,-e,-e)g\Big]+(-e,o,-e)h\Bigg\}$

+$(o,o,$-e*o*e$)\Bigg\{(o,o,e)\Big[(o,e,o)f+(-e,o,-e)g\Big]+(o,-e,-e)h\Bigg\}$

+$(o,o,N)\Bigg\{(o,o,e)\Big[(o,e,o)f+(-e,o,-e)g\Big]+(-e,o,-e)h\Bigg\}$

+$(o,o,N)\Big]\Bigg\{(o,o,e)\Big[(o,e,o)f+(-e,o,-e)g\Big]+(-e,-e,o)h\Bigg\}$

+$(o,o,N)\Bigg\{(o,o,e)\Big[(o,e,o)f+(-e,-e,o)g\Big]+(o,-e,-e)h\Bigg\}$

+$(o,o,N)\Bigg\{(o,o,e)\Big[(o,e,o)f+(-e,-e,o)g\Big]+(-e,o,-e)h\Bigg\}$

+$(o,o,$-e*o*e$)\Bigg\{(o,o,e)\Big[(o,e,o)f+(-e,-e,o)g\Big]+(-e,-e,o)h\Bigg\}$

+$(o,o,e)\Bigg\{(o,o,e)\Big[(o,o,e)f+(o,-e,-e)g\Big]+(-e,o,-e)h\Bigg\}$

+$(o,o,-e)\Bigg\{(o,o,e)\Big[(o,o,e)f+(o,-e,-e)g\Big]+(-e,-e,o)h\Bigg\}$

+$(o,o,N)\Bigg\{(o,o,e)\Big[(o,o,e)f+(-e,o,-e)g\Big]+(o,-e,-e)h\Bigg\}$

+$(o,o,$-e*o*e$)\Bigg\{(o,o,e)\Big[(o,o,e)f+(-e,o,-e)g\Big]+(-e,o,-e)h\Bigg\}$

+$(o,o,N)\Bigg\{(o,o,e)\Big[(o,o,e)f+(-e,o,-e)g\Big]+(-e,-e,o)h\Bigg\}$

+$(o,o,$-e*o*e$)\Bigg\{(o,o,e)\Big[(o,o,e)f+(-e,-e,o)g\Big]+(o,-e,-e)h\Bigg\}$

+$(o,o,N)\Bigg\{(o,o,e)\Big[(o,o,e)f+(-e,-e,o)g\Big]+(-e,o,-e)h\Bigg\}$

+$(o,o,N)\Big]\Bigg\{(o,o,e)\Big[(o,o,e)f+(-e,-e,o)g\Big]+(-e,-e,o)h\Bigg\}$

{\bf Definition~7.2} High Commutation Operators. Assume there is an
operator of two variables $\phi$, $ y_{1}\phi y_{2}=y_{0}$, its
commutation operators ~$\phi$(1,2,0),$\phi$(1,0,2),$\phi$(0,
\\2,1),$\phi$(2,1,0),$\phi$(2,0,1),$\phi$(0,1,2) will be defined by following formulas
and we introduce high commutation operators $\overline{C}$[1,2,0],
$\overline{C}$[1,0,2], $\overline{C}$[0,2,1], $\overline{C}$[2,1,0],
$\overline{C}$[2,0,1],$\overline{C}$[0,1,2] then new operators can
be expressed by $\phi$ and high commutation operators.
\begin{subequations}
\begin {equation}\label{eq:eps}
y_{1}\phi[1,2,0]y_{2}=y_{0} \qquad \qquad \phi[1,2,0]=
\overline{C}[1,2,0](\phi)
\end{equation}
\begin {equation}
y_{1}\phi[1,0,2]y_{0}=y_{2}\qquad\qquad\phi[1,0,2]=
\overline{C}[1,0,2](\phi)
\end{equation}
\begin {equation}
y_{0}\phi[0,2,1]y_{2}=y_{1}\qquad \qquad\phi[0,2,1]=
\overline{C}[0,2,1](\phi)
\end{equation}
\begin {equation}
y_{2}\phi[2,1,0] y_{1}=y_{0}\qquad \qquad\phi[2,1,0]=
\overline{C}[2,1,0](\phi)
\end{equation}
\begin {equation}
y_{2}\phi[2,0,1] y_{0}=y_{1}\qquad \qquad\phi[2,0,1]=
\overline{C}[2,0,1](\phi)
\end{equation}
\begin {equation}
y_{0}\phi[0,1,2]y_{1}=y_{2}\qquad \qquad\phi[0,1,2]=
\overline{C}[0,1,2](\phi)
\end{equation}
\end{subequations}

There are only two high commutation functions for a unary operator:
\begin{subequations}
\begin {equation}
\zeta_{e}(y)=y_{0}  \qquad \qquad\qquad   \zeta_{e}=\zeta
\end{equation}

\begin {equation}
\zeta_{t}(y_{0})=y   \qquad \qquad\qquad  \zeta_{t}=
\overline{C}(\zeta)
\end{equation}
\end{subequations}

{\bf Definition~7.3} High Tension-compression Operator.Assume there
is a binary operator $\phi$ and an unary operator $\zeta$, $\zeta
(y_{1})\phi y_{2}=y_{0}$, we can introduce a new binary operator
$\phi _{1}$ by $\phi$ and $\zeta$, $\phi _{1}$ will meet the
relation: $y_{1}\phi_{1} y_{2}=y_{0}$, that is say, $y_{1}\phi_{1}
y_{2}=\zeta(y_{1})\phi y_{2}$. Introduce a special high operator
$\overline{T}_{1}$ to express the relation between  $\phi _{1}$and
$\phi$,$\zeta$ .
\begin{subequations}
\begin {equation}
\phi_{1}=\phi \overline{T}_{1}\zeta
\end{equation}
In the same way if $y_{1}\phi \zeta(y_{2})=y_{0}$, we can introduce
a new binary operator $\phi _{2}$ by $\phi$ and $\zeta$, $\phi _{2}$
will meet the relation: $y_{1}\phi_{2} y_{2}=y_{0}$, that is say,
$y_{1}\phi_{2} y_{2}=y_{1}\phi \zeta (y_{2})$. Introduce a special
high operator $\overline{T}_{2}$ to express the relation between  $\phi
_{2}$and $\phi$,$\zeta$.
\begin {equation}
\phi_{2}=\phi \overline{T}_{2}\zeta
\end{equation}
If $y_{1}\phi y_{2}= \zeta (y_{0})$,that is say $\zeta
^{-1}[y_{1}\phi y_{2}]=y_{0}$, we can introduce a new binary
operator $\phi _{0}$ by $\phi$ and $\zeta$, $\phi _{0}$ will meet
the relation: $y_{1}\phi_{0} y_{2}=y_{0}$, that is say
$y_{1}\phi_{0} y_{2}=\zeta ^{-1}[y_{1}\phi y_{2}]$, and there is
$\overline{T}_{0}$:
\begin {equation}
\phi_{0}=\phi \overline{T}_{0}\zeta
\end{equation}
\end{subequations}

For an unary operator we have only $\overline{T}$ and
$\overline{T}_{0}$:
\begin{subequations}
\begin {equation}
\zeta_{1}\overline{T}\zeta_{2}=\zeta_{1}\zeta_{2}
\end{equation}
\begin {equation}
\zeta_{1}\overline{T}_{0}\zeta_{2}=\zeta_{2}^{-1}\zeta_{1}
\end{equation}
\end{subequations}

{\bf Definition~7.4}High Superposition Operator. Assume there are P operators of many
variables $\phi_{k}$ (k=1,p), its superposition operator $\phi$ will
be:
\begin {equation}
\phi=\sum_{k=1}^{P}\phi_{k}
\end{equation}

Function of $\phi$ will be the sum of function of $\phi_{k}$(k=1,p).
$\phi$ will be no-valued in a point if any of $\phi_{k}$  is
no-valued in this point. $\phi_{1}+\phi_{2}$ will be many-valued in
a point if $\phi_{1}$  is single-valued and $\phi_{2}$  is
many-valued in this point.

{\bf Definition~7.5}High Decomposition Operator.
\begin {equation}
\phi_{3}=\sum_{i=1}^{27}\zeta_{i}\Bigg\{\eta_{i4}\Big[\eta_{i1}(\widetilde{\beta_{1}})+\eta_{i2}(\widetilde{\beta_{2}})\Big]+\eta_{i3}(\widetilde{\beta_{3}})\Bigg\}
\end{equation}

 We can express the relations between$\zeta_{i}$ or
~$\eta_{ij}$ and~$\phi_{3}$with special operators $\overline{V}_{i}$
and $\overline{P}_{ij}$  and actually $\eta_{ij}$ is not change with
$\phi_{3}$.

\begin{subequations}
\begin {equation}
\zeta_{i}=\overline{V}_{i}(\phi_{3})  \qquad   \qquad  \qquad
\qquad (i=1, 27)
\end{equation}

\begin {equation}
\eta_{ij}=\overline{P}_{ij}(\phi_{3})        \qquad    \qquad (i=1,
27,j=1,4)
\end{equation}
\end{subequations}

Please note high commutation operator or high tension-compression
operator or high decomposition operator or high superposition
operator will be close within all three numbers operators if they
contain ones being many-valued and no-valued.

{\bf Definition~7.6}False operator of M+K variables. We  can change
an operator of M variables to a false one of (M+K ) variables by
adding $\sigma\widetilde{\zeta_{k}}$ in which $\sigma$ is a zero
operator and operator $\phi$  will not change with K variables.

\qquad\qquad\begin {equation} \phi=\sum_{i=1}^{L}fi\sum
_{j=1}^{M}g_{ij}\widetilde{\zeta_{i}}=\sum_{i=1}^{L}fi\Bigg\{\sum
_{j=1}^{M}g_{ij}\widetilde{\zeta_{i}}+\sum
_{k=M+1}^{M+K}\sigma\widetilde{\zeta_{k}}\Bigg\}
\end{equation}

We can also get false operator of (M+K ) variables from one of M
variables by $\overline{T}_{k}\sigma$  (k=M,M+K).

Formula solution of double branches operator equation with
parameterized operators will be:

$$y=\sum_{k=1}^{27}\overline{V}_{k}\Big[\overline{C}(2,3,0,1) \Big(\sum_{i=1}^{9}\big\{\big[\phi_{1}T_{3}o\overline{T}_{0}(\overline{P}_{i1}\phi_{3})^{-1}+\phi_{2}\overline{T}_{2}\sigma\overline{T}_{0}(\overline{P}_{i2}\phi_{3})^{-1}\big]\overline{T}_{0}(\overline{V}_{i}\phi_{3})^{-1}\big\}\Big)\Big]$$

$$\Bigg[\overline{P}_{k4}\Big[\overline{C}(2,3,0,1) \Big(\sum_{i=1}^{9}\big\{\big[\phi_{1}T_{3}o\overline{T}_{0}(\overline{P}_{i1}\phi_{3})^{-1}+\phi_{2}\overline{T}_{2}\sigma\overline{T}_{0}(\overline{P}_{i2}\phi_{3})^{-1}\big]\overline{T}_{0}(\overline{V}_{i}\phi_{3})^{-1}\big\}\Big)\Big]$$

$$\Bigg(\overline{P}_{k1}\Big[\overline{C}(2,3,0,1) \Big(\sum_{i=1}^{9}\big\{\big[\phi_{1}T_{3}o\overline{T}_{0}(\overline{P}_{i1}\phi_{3})^{-1}+\phi_{2}\overline{T}_{2}\sigma\overline{T}_{0}(\overline{P}_{i2}\phi_{3})^{-1}\big]\overline{T}_{0}(\overline{V}_{i}\phi_{3})^{-1}\big\}\Big)\Big](f)$$

$$+\overline{P}_{k2}\Big[\overline{C}(2,3,0,1) \Big(\sum_{i=1}^{9}\big\{\big[\phi_{1}T_{3}o\overline{T}_{0}(\overline{P}_{i1}\phi_{3})^{-1}+\phi_{2}\overline{T}_{2}\sigma\overline{T}_{0}(\overline{P}_{i2}\phi_{3})^{-1}\big]\overline{T}_{0}(\overline{V}_{i}\phi_{3})^{-1}\big\}\Big)\Big](g)\Bigg)$$

$$+\overline{P}_{k3}\Big[\overline{C}(2,3,0,1) \Big(\sum_{i=1}^{9}\big\{\big[\phi_{1}T_{3}o\overline{T}_{0}(\overline{P}_{i1}\phi_{3})^{-1}+\phi_{2}\overline{T}_{2}\sigma\overline{T}_{0}(\overline{P}_{i2}\phi_{3})^{-1}\big]\overline{T}_{0}(\overline{V}_{i}\phi_{3})^{-1}\big\}\Big)\Big](h)\Bigg]$$

Please note solution for double branches operator equation has the
same form with one for double branches algebraic equation. Is it
appropriate to class algebraic equation and operator equation to
different fields? But we have done it! Mathematics has been parted
to many alone islands. This situation is not good and will be
changed in future.

These results mean that there is a new accurate analytical route
beside approximate numerical method and topological way in study of
operator equations certainly including functional equations and
function equations and differential equations.We can extend results
to N numbers operators of M variables but there are many works to
be done.

\section{Try to extend to continuous situation}

    We can extend results about discrete functions to continue
functions if we accept results about Hilbert's 13th problem. We can
express formula solution of equation constructed by continue
functions in the same form of equation constructed by discrete
functions even though we can't give a procedure to decompose a
continue function of many variables to a superposition of functions
of one variable. But there are many tasks to be done if we want to
make results to be strict in logic.

We must prove that every continue operator of many variables can be
represented as a superposition of continue operators of one variable
if we want to extend results in this paper to continue operators and
equations constructed by them. I don't know if there is such a
result in current literature. Please give it if there isn't.

 There are enough space for us to write our results so we are
luckier than Pierre de Fermat (1601-1665) who could not write
the proof of his last theorem. Now we have only poor results shown here
but mathematicians will find more and more good results because
there is huge mineral deposit in this direction. Please believe this
point!



\end{document}